\documentclass[11pt]{article}
\usepackage{amsthm}
\usepackage[english]{babel}
\usepackage[latin1]{inputenc}

\usepackage{amsmath,amssymb}

\usepackage{epsf}
\usepackage[all]{xy}

\usepackage{amsmath,amsthm,amssymb,amsfonts,indentfirst}

\newcommand{\rmap}       {\longrightarrow}

\newcommand{\frakg}     {\mathfrak{g}}

\theoremstyle{plain}
\newtheorem{theorem}{Theorem}[section]
\newtheorem{proposition}{Proposition}[section]
\newtheorem{lemma}{Lemma}[section]
\newtheorem{corollary}{Corollary}[section]
\theoremstyle{definition}

\DeclareMathAlphabet{\mathpzc}{OT1}{pcz}{m}{it}

\begin{document}

\title{Multiplicative Dirac structures on Lie groups}

\author{Cristi\'an Ortiz\footnote{This research was supported by Capes-Brazil via a PEC-PG scholarship.}\\ \small{\textit{Instituto Nacional de Matem\'atica Pura e Aplicada}}\\ \small{\textit{cortiz@impa.br}}}
\date{}

\maketitle
\begin{abstract}
\indent We study multiplicative Dirac structures on Lie groups. We show that the characteristic foliation of a multiplicative Dirac structure is given by the cosets of a normal Lie subgroup and, whenever this subgroup is closed, the leaf space inherits the structure of a Poisson-Lie group. We also describe multiplicative Dirac structures on Lie groups infinitesimally. 

\vskip 0.5\baselineskip

\hspace{6.1cm}{\bf R\'esum\'e}
\vskip 0.5\baselineskip
\noindent
Nous \'etudions les structures de Dirac multiplicatives sur les groupes de Lie. On montre que le feuilletage caract\'eristique d'une structure de Dirac multiplicative est donn\'ee par les classes \`a gauche (respectivement \`a droite) d'un sous-groupe distingu\'e et, quand ce sous-groupe est ferm\'e, l'espace des feuilles est muni d'une structure de groupe de Lie-Poisson. Nous d\'ecrivons aussi la version infinit\'esimale des structures de Dirac multiplicatives sur les groupes de Lie.  

\end{abstract}

\section*{Version fran\c{c}aise abr\'eg\'ee}
Un groupe de Lie-Poisson est un objet du type groupe de Lie dans la cat\'egorie des vari\'et\'es de Poisson, c'est-\`a-dire un groupe de Lie $G$ muni d'une structure de Poisson $\pi\in \Gamma(\bigwedge^{2}TG)$ telle que la multiplication $m:G\times G \longrightarrow G$ soit une application de Poisson, le produit $G\times G$ \'etant muni de la structure de Poisson produit. De mani\`ere \'equivalente, le tenseur de Poisson satisfait la condition multiplicative $\pi_{gh}=(L_g)_{*}\pi_h + (R_h)_{*}\pi_g$. Les groupes de Lie-Poisson ont \'et\'e introduits par Drinfeld (\cite{D}). Ces structures apparaissent aussi dans l'etude des propri\'et\'es Hamiltonnienes des transformations d'habillage de certain syst\`emes int\'egrables(\cite{STS}).

La donn\'ee infinit\'esimale associ\'ee \`a un groupe de Lie-Poisson $G$ est son alg\`ebre de Lie $\mathfrak{g}$ accompagn\'ee d'une structure suppl\'ementaire: l'espace dual $\mathfrak{g}^{*}$ h\'erite d'une structure d'alg\`ebre de Lie, satisfaisant une condition de compatibilit\'e avec le crochet de Lie sur $\mathfrak{g}$ (cf. par exemple \cite{D,LW}). Une telle paire d'alg\`ebres de Lie $(\mathfrak{g},\mathfrak{g}^{*})$ s'appelle une bialg\`ebre de Lie. R\'eciproquement, toute bialg\`ebre de Lie $(\mathfrak{g},\mathfrak{g}^{*})$ est la bialg\`ebre de Lie d'un groupe de Lie-Poisson $G$ connexe et simplement connexe (\cite{D}). 

Une structure de Poisson peut \^etre vue comme un cas particulier d'une structure g\'eom\'etrique plus g\'en\'erale: la structure de Dirac. La notion de structure de Dirac a \'et\'e introduite par Courant et Weinstein dans le cadre des syst\`emes m\'echaniques (\cite{C,CW}). Les structures de Dirac incluent les structures pr\'e-symplectiques et les structures de Poisson, de m\^eme que les feuilletages r\'eguliers. Si $(G,\pi)$ est un groupe de Lie-Poisson, la multiplicativit\'e de $\pi$ est \'equivalente au fait que l'application fibr\'ee induite $\pi^{\sharp}: T^{*}G \longrightarrow TG$ soit un morphisme de groupo\"{i}des. Dans ces conditions, la structure de Dirac associ\'ee $L_{\pi}=$graph$(\pi^{\sharp})\subseteq TG\oplus T^{*}G$ est un sous-groupo\"{i}de de Lie.

Dans cet expos\'e, nous \'etudions les groupes de Lie $G$ \'equip\'es d'une structure de Dirac $L$ multiplicative  dans le sens o\`u $L\subseteq TG\oplus T^{*}G$ est un sous-groupo\"{i}de. Nous appelerons un tel couple $(G,L)$ un groupe de Lie-Dirac.

Cet expos\'e est organis\'e de la mani\`ere suivante: dans le deuxi\`eme paragraphe, nous \'etudions les feuilletages multiplicatifs sur un groupe de Lie $G$, c'est-\`a-dire les structures de Dirac multiplicatives induites par des distributions r\'eguli\`eres involutives sur $TG$. \'Etant donn\'e un sous-groupe de Lie connexe $H$ d'un groupe de Lie $G$, il existe un feuilletage canonique de $G$, dont les feuilles sont les classes \`a droite de $H$. Nous montrons que ce feuilletage est multiplicatif si et seulement si $H$ est un sous-groupe normal de $G$.

Dans le troisi\`eme paragraphe, nous \'etudions le feuilletage caract\'eristique d'une structure de Dirac multiplicative. Nous montrons que le noyau d'une structure de Dirac multiplicative sur un groupe de Lie $G$ est une distribution r\'eguli\`ere involutive, que s'int\`egre comme un feuilletage r\'egulier et multilicatif de $G$. \`A l'aide de la description des feuilletages multiplicatifs donn\'es au paragraphe $2$, nous d\'eduisons que le feuilletage caract\'eristique d'un groupe de Lie-Dirac $(G,L)$ est simple et que l'espace des feuilles est un groupe de Lie-Poisson. C'est-\`a-dire que les groupes de Lie-Dirac, apr\`es un quotient canonique deviennent des groupes de Lie-Poisson. Enfin, nous appliquons ce r\'esultat ainsi que la correspondance de Drinfeld entre les groupes de Lie-Poisson et les bialg\`ebres de Lie pour obtenir la version infinit\'esimale des groupes de Lie-Dirac.

\section{Introduction}
A Poisson-Lie group is a Lie group $G$ with a multiplicative Poisson structure $\pi\in \Gamma(\Lambda^{2}TG)$, that is, the multiplication map $m:G\times G\longrightarrow G$ is a Poisson map \cite{D,LW,STS}. Equivalently, the bundle map $\pi^{\sharp}: T^{*}G\longrightarrow TG$ is a groupoid morphism, where the tangent and cotangent bundles have the groupoid structures induced by $G$  (see e.g. \cite{M}). Therefore, the multiplicativity property of the Poisson bivector $\pi$ is equivalent to saying that the associated Dirac structure \cite{C,CW} $L_{\pi}=\mathrm{graph}(\pi^{\sharp})$ is a Lie subgroupoid of the direct-sum $\mathcal{VB}$-groupoid $\mathbb{T}G=TG\oplus T^{*}G$. This motivates the following definition: a Dirac structure $L$ on a Lie group $G$ is \textbf{multiplicative} if $L\subseteq \mathbb{T}G$ is a subgroupoid. We refer to a Lie group equipped with a multiplicative Dirac structure as a Dirac-Lie group.

Clearly Poisson-Lie groups are particular examples of Dirac-Lie groups. On the other extreme, one can check that there are no interesting multiplicative $2$-forms on Lie groups -- the only one is the zero $2$-form. Another class of examples is obtained as follows: Let $p:G_1\longrightarrow G_2$ be a homomorphism of Lie groups which is a surjective submersion. If $\pi$ is a multiplicative Poisson structure on $G_2$, then its pull back (in the sense of Dirac structures \cite{BR}) turns out to be a multiplicative Dirac structure on $G_1$, whose pre-symplectic leaves are the inverse images by $p$ of the symplectic leaves of $G_2$, and whose characteristic foliation is given by the fibres of the submersion $p$. Our main observation in this note is that, modulo a regularity condition, all multiplicative Dirac structures on Lie groups are of this form.




\section{Preliminaries}
\subsection{Multiplicative foliations}

Let $G$ be a Lie group. If $F\subseteq TG$ is a regular integrable distribution, one can check that the corresponding Dirac structure $L_F=F\oplus \mathrm{Ann}(F)$ is multiplicative if and only if $F\subseteq TG$ is a Lie subgroup, where $TG$ has the natural Lie group structure induced from $G$.\\

\begin{proposition}\label{mf}
 Let $\mathcal{F}$ be the foliation integrating a multiplicative distribution $F\subseteq TG$. The following holds:

\begin{enumerate}
 \item The leaf through the identity $\mathcal{F}_e\subseteq G$ is a normal Lie subgroup.

 \item The foliation $\mathcal{F}$ is given by cosets of $\mathcal{F}_e$.\\

\end{enumerate}
\end{proposition}

\noindent\textbf{Proof.} Since $F \subseteq TG$ is a subgroup, it is closed under
multiplication in $TG$, that is $\rm{d}m(g, h)(X_g,X_h) = \rm{d}R_{h}(g)X_g
+\rm{d}L_{g}(h)X_h \in F_{gh}$ for every $X_g,X_h \in F$. In
particular, for $X_h = 0$ we see that $F$ is right invariant, i.e.
$\rm{d}R_{h}(g)X_g \in F_{gh}$. Similarly we obtain left invariance of
$F$: $\rm{d}L_{g}(h)X_h \in F_{gh}$. This says that the distribution at each $g\in G$ is given by

\begin{equation}\label{eq:1}
 F_g = \rm{d}L_{g}(e)F_e = \rm{d}R_{g}(e)F_e.
\end{equation}

Consider now $\mathcal{F}_e$, the leaf of $F$ through the identity $e \in
G$. For every $a, b\in \mathcal{F}_e$ there exist paths $a(t), b(t) \in G, t\in [0,1]$, tangent to the
distribution $F$, joining the identity $e\in G$ to $a$ and $b$,
respectively. We want to prove that $c = ab\in \mathcal{F}_e$. For this,
take the path $c(t) = a(t)b(t)$, which joins the identity to $c =
ab$. The path $c(t)$ is
tangent to the distribution $F$: indeed, the bi-invariance of $F$ implies that
$$
 c'(t) =\rm{d}R_{b}(t)(a(t))a'(t) + \rm{d}L_{a}(t)(b(t))b'(t) \in F_{c(t)},
$$
since $a'(t) \in F_{a(t)}$ and $b'(t) \in F_{b(t)}$. This shows that $c\in\mathcal{F}_e$. A
similar computation shows that $\mathcal{F}_{e}$ is closed by the inversion map. Therefore the leaf through the identity is a subgroup of $G$. Moreover, it follows from (\ref{eq:1}) that the Lie algebra of $\mathcal{F}_e$ is Ad-invariant, which is equivalent to $\mathcal{F}_e$ being a normal subgroup. The assertion in $(ii)$ follows from the bi-invariance in (\ref{eq:1}).\qed


\subsection{Functorial properties of multiplicative Dirac structures}

Let $M$ be a smooth manifold. The generalized tangent bundle of $M$ is the direct-sum vector bundle $\mathbb{T}M=TM\oplus T^{*}M$. Given a smooth map $\varphi:M_1\rmap M_2$ and $x\in M_1$, we say that the elements $\eta=(X,\alpha)\in (\mathbb{T}M_1)_x$ and $\xi=(Y,\beta)\in (\mathbb{T}M_2)_{\varphi(x)}$ are \textbf{$\varphi$-related} if $Y=\rm{d}\varphi(X)$ and $\alpha=\rm{d}\varphi^{*}\beta$. For a Lie group $G$ with Lie algebra $\frakg$, $TG$ is a Lie group, $T^{*}G$ is a Lie groupoid over $\frakg^{*}$ and we consider the groupoid $\mathbb{T}G=TG\oplus T^{*}G$; we denote the groupoid multiplication in $\mathbb{T}G$ by $\xi*\xi'$, when $\xi,\xi'$ are composable.\\

\begin{lemma}\label{lem:1}
Let $\varphi: G_1\rmap G_2$ be a homomorphism of Lie groups, which
is a surjective submersion. If $\xi_{\varphi(g)}, \xi'_{\varphi(h)}\in \mathbb{T}G_2$ are $\varphi$-related to $\eta_{g},\eta'_{h}\in \mathbb{T}G_1$, respectively, then $\xi_{\varphi(g)},\xi'_{\varphi(h)}$ are composable if and only if $\eta_g, \eta'_h$ are composable. Moreover,  $\xi_{\varphi(g)}*\xi'_{\varphi(h)}\in \mathbb{T}G_2$ is $\varphi$-related to $\eta_{g}*\eta'_{h}\in \mathbb{T}G_1$.\\

\end{lemma}


\textbf{Proof.} The pull back bundles $\varphi^{*}(TG_2)$, $\varphi^{*}(T^{*}G_2)$ have natural groupoid structures in such a way that $\rm{d}\varphi:TG_1\rmap \varphi^{*}(TG_2)$ and $\rm{d}\varphi^{*}:\varphi^{*}(T^{*}G_2)\rmap T^{*}G_1$ are morphisms of groupoids. The statements follow from this fact and a direct computation using that $\varphi$ is a surjective submersion. \qed\\

\begin{corollary}\label{cor:func}
 Let $\varphi: G_1\rmap G_2$ be a homomorphism of Lie groups, which
is a surjective submersion. Assume that $L_1,L_2$ are Dirac
structures on $G_1,G_2$, respectively. If $\varphi$ is a forward Dirac map and $L_1$ is multiplicative, then $L_2$ is multiplicative. Also, if $\varphi$ is a backward Dirac map and $L_2$ is multiplicative, then $L_1$ is multiplicative.\\ 
\end{corollary}

 
\textbf{Proof.} Recall that $\varphi$ is a forward Dirac map if and only if $L_2$ is the bundle of all $\varphi$-related elements to elements in $L_1$. A backward Dirac map is defined in a similar way, see \cite{BR}. The statement follows from Lemma \ref{lem:1}.\qed


\section{The main result}

Let $L$ be a Dirac structure on a smooth manifold $M$. Let us denote by $p_{T}$ the canonical projection of $TM\oplus T^{*}M$ on $TM$. The manifold $M$ carries a singular foliation tangent to the generalized distribution $p_{T}(L)\subseteq TM$, and each leaf of this foliation inherits a canonical pre-symplectic structure whose kernel is ker$(L)=L\cap TM$. If the distribution ker$(L)$  has constant rank, then it is integrable. If the corresponding foliation $\mathcal{K}$ is simple, then the leaf space $M/\mathcal{K}$ has a natural Poisson structure obtained by identifying functions on the quotient with \textbf{admissible functions} on $M$, see \cite{C} for details. The foliation $\mathcal{K}$ is called the \textbf{characteristic foliation} of $L$. 

In the special case of Dirac-Lie groups our main result is the following.\\

\begin{theorem}\label{thm:1}
Let $G$ be a Lie group with a multiplicative Dirac structure
$L\subseteq TG\oplus T^{*}G$. Then:

\begin{enumerate}
 \item The kernel of $L$ is a multiplicative integrable distribution, and the leaves of the characteristic foliation $\mathcal{K}$ are cosets of the normal Lie subgroup $\mathcal{K}_e\subseteq G$.

 \item If $\mathcal{K}_e$ is closed\footnote{Notice that not always the characteristic leaf through the identity $\mathcal{K}_e\subseteq G$ is a closed subgroup. Indeed, if $G$ is the torus and $K\subseteq G$ is a dense geodesic, then the Dirac structure associated to the foliation of $G$ by cosets of $K$ is a multiplicative Dirac structure whose characteristic leaf through the identity is dense.}, then the leaf space $G/\mathcal{K}$ is smooth and
the induced Poisson structure $\pi$  is multiplicative (i.e., $G/\mathcal{K}$ becomes a Poisson-Lie group).
Moreover, $L$ is the pull back of $\pi$ by the quotient map $G\rmap G/\mathcal{K}$.\\
\end{enumerate}

\end{theorem}

 
\textbf{Proof.} Since $L$ is multiplicative, we have that $\ker(L)=L\cap TG\subseteq
TG$ is a subgroup, hence (\ref{eq:1}) implies that $\ker(L)$ has constant rank. In particular it defines an involutive distribution, whose leaves are given by cosets of the normal Lie subgroup $K=\mathcal{K}_e$ (the leaf through the identity) by Prop.~\ref{mf}. If $K$ is closed, then $G/K$ is a Lie group and the projection $G\rmap G/K$ is a surjective submersion which is both a forward and backward Dirac map \cite{BR}, where $G/K$ is equipped with the natural Poisson structure induced by $L$. The multiplicativity property of this Poisson structure is a consequence of Corollary \ref{cor:func}.\qed\\


Combining Theorem \ref{thm:1} and Drinfeld's correspondence between Poisson-Lie groups and Lie bialgebras \cite{D}, we obtain the infinitesimal counterpart of Dirac-Lie groups.\\

\begin{corollary}\label{cor:gi}

 Let $G$ be a Lie group with Lie algebra $\frakg$. If $G$ is equipped with a multiplicative Dirac structure $L$, then $\mathfrak{k}=$ker$(L)_e$ is an ideal in $\mathfrak{g}$ and the quotient $\frakg/\mathfrak{k}$ inherits the structure of a Lie bialgebra.\\

\end{corollary}

\textbf{Proof.} The multiplicativity of the characteristic distribution implies that $\mathfrak{k}\subseteq \frakg$ is an ideal. Now consider the connected and simply connected Lie group $T$ integrating the quotient Lie algebra $\frakg/\mathfrak{k}$. The canonical projection $\frakg \longrightarrow \frakg/\mathfrak{k}$ integrates to a homomorphism of Lie groups $\phi: \tilde{G}\rmap T$, where $\tilde{G}$ denotes the universal covering of $G$. The subgroup $H=\ker(\phi)$ is closed and normal in $\tilde{G}$, therefore the connected component of the identity $H_0$ is closed and normal as well and the quotient group $\tilde{G}/H_0$ inherits a Poisson-Lie structure. Since $\tilde{G}/H$ is locally diffeomorphic to $\tilde{G}/H_0$, the Lie algebra $\frakg/\mathfrak{k}$ inherits a Lie bialgebra structure.\qed\\
In the situation of Corollary \ref{cor:gi}  we say that $(G,L)$ is an \textbf{integration} of the infinitesimal data $(\frakg,\mathfrak{k})$, where $\mathfrak{k}\subseteq \frakg$ is ideal and $\frakg/\mathfrak{k}$ is a Lie bialgebra.\\

\begin{corollary}
 If $G$ is connected and simply connected and $\mathfrak{k}\subseteq \frakg$ is an ideal such that $\frakg/\mathfrak{k}$ is a Lie bialgebra, then there is a unique multiplicative Dirac structure on $G$ integrating $(\frakg,\mathfrak{k})$.\\
\end{corollary}

 \textbf{Proof.} Let $T$ be the connected and simply connected Lie group integrating $\frakg/\mathfrak{k}$. Consider the homomorphism $\phi:G\rmap T$ and $H\subseteq G$ as in the proof of Corollary \ref{cor:gi}. The quotient group $G/H\cong T$ has a multiplicative Poisson structure $\pi_T$ integrating the Lie bialgebra $\frakg/\mathfrak{k}$. Since $\phi$ is a surjective submersion, we induce a multiplicative Dirac structure $L$ on $G$ according to Corollary \ref{cor:func}. This shows that $(G,L)$ is an integration of $(\frakg,\mathfrak{k})$.\qed 

\section*{Acknowledgements}
I would like to thank Henrique Bursztyn for helpful discussion and valuable comments. I thank also Marius Crainic for helpful suggestions. Finally I thank the Department of Mathematics at Utrecht University for hospitality while part of this work was carried out.

\end{document}